\documentclass[12pt]{amsart}
\usepackage[utf8]{inputenc}
\usepackage[T2A,T1]{fontenc}
\usepackage[russian,english]{babel}
\usepackage{textcase}

\usepackage[all, cmtip]{xy}
\usepackage{amsfonts,amsmath,amstext,amsbsy,amssymb, amsopn,amsthm,upref,eucal, mathrsfs}

\newtheorem{proposition}{Proposition}
\newtheorem{theorem}{Theorem}

\newtheorem{corollary}{Corollary}
\newtheorem{remark}{Remark}

\DeclareMathOperator{\Conf}{Conf}
\DeclareMathOperator{\supp}{supp}

\begin{document}

\title[Stochastic Euler Products]{Stochastic Euler Products in the Paley--Wiener Space}

\author{Alexander\ I.~Bufetov}
\address{Steklov Mathematical Institute of Russian Academy of Sciences, Moscow}
\email{bufetov@mi-ras.ru}

\thanks{This work was supported by the Russian Science Foundation (grant No. 24-71-10109).}

\date{}

\begin{abstract}
General sufficient conditions, in terms of stochastic Euler products, are given for a realization of a point process to have finite excess or finite deficiency in the Paley-Wiener space. 
\end{abstract}

\maketitle

\begin{center}
  
\it{
To Andrei Yurevich Zaitsev, with profound respect}

\end{center}

\section{Introduction}

Almost every realization of a determinantal point process associated with an orthogonal projection that is locally of trace class is a uniqueness set for the range of the projection \cite{BQS}. Equivalently, the family of reproducing kernels of the space defining our determinantal process, evaluated at the points of a realization of the process, is almost surely complete. In the special case of the sine process, defined by the projection onto the Paley--Wiener space $PW=PW([-\pi, \pi])$ of square-integrable functions whose Fourier transforms are supported in the interval $[-\pi,\pi]$, the Completeness Theorem is due to Ghosh \cite{ghosh-exp}; in the case of discrete phase spaces, it is due to Lyons \cite{lyons-det}. The proof of completeness given in \cite{BQS} relies essentially on the determinantal structure.

The Completeness Theorem from \cite{BQS} raises the question of the minimality of our point sets---in other words, whether completeness is preserved when a particle is removed from our random configuration. Minimality clearly holds in the finite-dimensional case. However, for processes with infinitely many particles, minimality fails: the sine process almost surely has excess one \cite{bufetov-excess}. Almost every realization of the sine process remains complete after one particle is removed, but after two particles are removed the resulting point set is the zero set of a nonzero function in the Paley--Wiener space. The second statement---that the excess is at most one---was established in \cite{bufetov-excess} for all determinantal processes with integrable kernels $\Pi$ such that
\begin{equation}
\int\limits_{\mathbb R} \frac{\Pi(x,x)}{1+x^2}\,dx<+\infty.
\end{equation}

The aim of this note is to give general sufficient conditions ensuring that a realization of a point process belongs to the Paley--Wiener space; that it belongs to this space after a finite (deterministic!) number of particles are removed; and, conversely, sufficient conditions ensuring that a realization does not belong to the Paley--Wiener space.

Informally, the results of the note can be summarized as follows. Suppose initially that our stochastic Euler product grows at infinity no faster than $\exp(\pi z)$. If a nontrivial limit theorem holds for the logarithm of the absolute value of the stochastic Euler product, then, with positive probability, realizations of our point process are uniqueness sets for the Paley--Wiener space. If, on the other hand, the Euler product itself converges to a random distribution, then, with positive probability, realizations of our point process have excess at least one. While this note is devoted to the Paley-Wiener space only, its methods are quite general, and should be applicable to more general integral transforms, such as, for instance,  the Hankel (cf., e.g., \cite{Akh54}), the  Airy, the Kontorovich-Lebedev transform \cite{KL38} as well as the Gorbunov transform in the theory of confluent hypergeometric functions \cite{Gor25}.

We now turn to the precise statements. Let $\Conf(\mathbb{R})$ denote the space of configurations on $\mathbb{R}$. Recall that a configuration $X\in\Conf(\mathbb{R})$ is called a uniqueness set, or equivalently a complete set, for the Paley--Wiener space if a function in the Paley--Wiener space is uniquely determined by its restriction to $X$---or, equivalently, if a function in the Paley--Wiener space that vanishes on $X$ is identically zero. A uniqueness set is minimal if it ceases to be a uniqueness set once a particle is removed. A uniqueness set has excess $\ell$ if, after removing $\ell$ distinct particles, it remains complete and becomes minimal. A configuration $X\in\Conf(\mathbb{R})$ has deficiency $\ell$ if adding $\ell$ distinct particles makes it a complete minimal set. Finally, we say that a configuration $X$ has finite deviation if it has either finite excess or finite deficiency.

Let $\mathscr{C}(\mathbb{R})$ denote the set of configurations $X\in\Conf(\mathbb{R})$ such that
\begin{enumerate}
    \item $\displaystyle\sum_{x\in X,\,x\neq 0}\frac{1}{x}$ is defined in the principal-value sense;
    \item $\displaystyle\sum_{x\in X,\,x\neq 0}\frac{1}{x^2}<+\infty$.
\end{enumerate}
Every complete minimal set in the Paley--Wiener space belongs to $\mathscr{C}(\mathbb{R})$ (Theorem 6, p.~544, in Levin \cite{Levin}).

For configurations $X\in\mathscr{C}(\mathbb{R})$ and numbers $z\in\mathbb{C}$, $w\in\mathbb{C}\setminus\mathbb{R}$, define the Euler product
\[
\widetilde{G}_X(z;w)=(\text{p.v.})\prod_{x\in X}\left(\frac{x-z}{x-w}\right),
\]
and set $G_X(z)=\widetilde{G}_X(z;\,i)$.

If there is an $a>0$ such that
\begin{equation}\label{type-pi-poly}
\sup\limits_{z\in\mathbb C}|G_X(z)|\exp(-\pi|z|)(1+|z|)^{-a}<+\infty,
\end{equation}
and the restriction of the Euler product $G_X$ to the real axis is square-integrable, then $X$ is not complete in the Paley--Wiener space.

\begin{remark}
The Paley--Wiener Theorem in the form of Plancherel and Pólya(cf. Theorem II, p.~226, in \cite{pl-po}) gives a slightly stronger statement: if $G_X$ is an entire function of exponential type at most $\pi$, that is, if
\begin{equation}\label{type-pi-expo}
\limsup\limits_{r\to\infty}\frac{\log\max\limits_{z:|z|=r}|G_X(z)|}{r}\leq\pi,
\end{equation}
and the restriction of the Euler product $G_X$ to the real axis is square-integrable, then $X$ is not complete in the Paley--Wiener space. Nevertheless, the Euler products of realizations of all point processes considered in this paper almost surely satisfy the more restrictive condition \eqref{type-pi-poly} as well.
\end{remark}

The converse is not true in general: the fact that the Euler product $G_X$ is not square-integrable does not by itself imply that the configuration $X$ is complete in the Paley--Wiener space. The set of even integers provides a counterexample. Nonetheless, if the Euler product $G_X$ has regular behaviour on the imaginary axis, then its excess or deficiency can indeed be calculated in terms of $G_X$.

\begin{proposition}\label{low-imag-prop}
Let $X\in\mathscr{C}(\mathbb{R})$ be a configuration such that its Euler product $G_X$ satisfies inequality \eqref{type-pi-poly}.

Suppose, in addition, that for some $A>0$, the Euler product $G_X$ satisfies
\begin{equation}\label{low-imag}
\liminf\limits_{v\in\mathbb R:\,|v|\to\infty}|G_X(iv)|\exp(-\pi|v|)(1+v^2)^A>0.
\end{equation}
Then the following statements hold:
\begin{enumerate}
\item If
\[
G_X(t)\notin L_2(\mathbb R)\qquad\text{and}\qquad\frac{G_X(t)}{t+i}\in L_2(\mathbb R),
\]
then $X$ is complete and minimal in the Paley--Wiener space.

\item If $\ell\in\mathbb N$ is such that
\[
\frac{G_X(t)}{(t+i)^\ell}\notin L_2(\mathbb R)\qquad\text{and}\qquad\frac{G_X(t)}{(t+i)^{\ell+1}}\in L_2(\mathbb R),
\]
then $X$ has excess $\ell$ in the Paley--Wiener space.

\item If $\ell\in\mathbb N$ is such that
\[
G_X(t)(t+i)^{\ell-1}\in L_2(\mathbb R)\qquad\text{and}\qquad G_X(t)(t+i)^\ell\notin L_2(\mathbb R),
\]
then $X$ has deficiency $\ell$ in the Paley--Wiener space.
\end{enumerate}
\end{proposition}

\begin{proof}
It suffices to prove the first statement; the second and the third follow. By the Paley--Wiener Theorem, if $p\in X$, then
\[
\frac{G_X(z)}{z-p}\in PW,
\]
so the set $X\setminus\{p\}$ is not complete in the Paley--Wiener space. Now suppose, to the contrary, that there exists a nonzero function $g\in PW$ such that $g|_X=0$. We follow B.Ya.\,Levin, cf.\,\cite[Proposition\,(\foreignlanguage{russian}{г}), p. 542]{Levin}. Set $h=g/G_X$. By definition, $h$ is of zero exponential type. In view of \eqref{low-imag}, $h$ has polynomial growth on the imaginary axis, and hence the Phragmén--Lindelöf principle \cite{phrag,phrag-linde} implies that $h$ is a polynomial. But the Euler product $G_X$ itself is not square-integrable on the real line, which is a contradiction. The Proposition is proved.
\end{proof}

Observe that, in order to define the absolute value $|G_X(iv)|$ of the Euler product on the imaginary axis, we need only the convergence of the sum of the inverse squares of the particles---that is, only the second condition in the definition of $\mathscr{C}(\mathbb{R})$. Moreover, the infimum in \eqref{low-imag} may be taken over integer points only:
\begin{proposition}\label{rtozinf}
For every configuration $X\in\Conf(\mathbb R)$ such that
\[
\displaystyle\sum_{x\in X,\,x\neq 0}\frac{1}{x^2}<+\infty,
\]
and every $A\in\mathbb R$, we have
\begin{multline}\label{low-imag-z}
\liminf\limits_{v\in\mathbb R:\,|v|\to\infty}|G_X(iv)|\exp(-\pi|v|)(1+v^2)^A\geq\\
\geq\exp(-\pi)\liminf\limits_{k\in\mathbb Z:\,|k|\to\infty}|G_X(ik)|\exp(-\pi|k|)(1+k^2)^A.
\end{multline}
\end{proposition}
This is clear.

Now let $\mathbb{P}$ be a point process such that $\mathbb{P}(\mathscr{C}(\mathbb{R}))=1$. The random entire function $G_X$ is called the stochastic Euler product associated with the point process $\mathbb{P}$. Our first aim is to give simple sufficient conditions under which the stochastic Euler product almost surely enters the Paley--Wiener space after a finite (deterministic!) number of particles are removed and, conversely, leaves the Paley--Wiener space after a finite (deterministic!) number of particles are added.

\begin{proposition}\label{small}
Suppose that, for some $N_0>0$ and $B>0$, the following estimates hold:
\begin{equation}\label{rad-gr}
\limsup\limits_{m\to+\infty}\exp(-2\pi m)\cdot(m^2+1)^{-N_0}\mathbb{E}_{\mathbb{P}}\int\limits_{|z|=m}|G_X(z)|^2\frac{dz}{iz}<+\infty;
\end{equation}
\begin{equation}\label{sm-l2}
\mathbb{E}_{\mathbb{P}}\int\limits_{\mathbb R}\frac{|G_X(u)|^2}{(u^2+1)^B}\,du<+\infty.
\end{equation}
Then, for every $\ell\in\mathbb N$ with $\ell\geq B$, for $\mathbb{P}$-almost every configuration $X$ and every choice of $\ell$ distinct particles $p_1,\ldots,p_\ell\in X$, we have
\[
\frac{G_X(z)}{(z-p_1)\cdots(z-p_\ell)}\in PW[-\pi;\,\pi].
\]
In particular, for $\mathbb{P}$-almost every configuration $X$, either $X$ is itself incomplete in $PW$, or it has excess at most $\ell-1$ in $PW$.
\end{proposition}

Proposition \ref{small} directly implies
\begin{corollary}\label{small-ptws}
Suppose that estimate \eqref{rad-gr} holds and, in addition, that for some $A>0$ and every $u\in\mathbb R$,
\[
\mathbb{E}_{\mathbb{P}}|G_X(u)|^2\leq(u^2+1)^A.
\]
Then, for every $\ell\in\mathbb N$ with $\ell>A+1/2$, for $\mathbb{P}$-almost every configuration $X$ and every choice of $\ell$ distinct particles $p_1,\ldots,p_\ell\in X$, we have
\[
\frac{G_X(z)}{(z-p_1)\cdots(z-p_\ell)}\in PW[-\pi;\,\pi].
\]
In particular, for $\mathbb{P}$-almost every configuration $X$, either $X$ is itself incomplete in $PW$, or it has excess at most $\ell-1$ in the space $PW=PW[-\pi;\,\pi]$.
\end{corollary}

To obtain an upper bound on the deficiency of a random configuration $X$, we will need stronger assumptions.

We begin with a simple general observation.
\begin{proposition}\label{bdd-dev}
If $X\in\mathscr{C}(\mathbb{R})$ has finite deviation, then $X$ is a uniqueness set in the Paley--Wiener space if and only if $G_X\notin L_2(\mathbb R)$.
\end{proposition}
Indeed, the statement directly follows from Theorem 6, p.~544, in Levin \cite{Levin}.

We say that a point process $\mathbb{P}$ has bounded deviation if there exists an $N\in\mathbb N$ such that, for $\mathbb{P}$-almost every configuration, either its excess is at most $N$ or its deficiency is at most $N$.

\begin{proposition}\label{big}
Suppose that estimate \eqref{rad-gr} holds and, in addition, that for some $r_0>0$, $A>0$, $C>0$, and $\delta>0$, and for every $z\in\mathbb C$, $z=u+iv$, with $|v|\geq r_0$, we have
\begin{multline}\label{log-dev}
\mathbb{P}\!\left(\left\{X\in\mathscr{C}(\mathbb{R}):\left|\log|G_X(z)|-\pi|v|\right|\geq A\log(2+|z|)\right\}\right)\\
\leq\frac{C}{(1+|z|)^{1+\delta}}.
\end{multline}
Then the point process $\mathbb{P}$ has bounded deviation.
\end{proposition}

Observe that, in order to obtain estimate \eqref{log-dev}, it is sufficient to have the estimate
\[
\mathbb{E}_{\mathbb{P}}\!\left(|G_X(z)|^2\exp(-2\pi|v|)+|G_X(z)|^{-2}\exp(2\pi|v|)\right)\leq C(u^2+1)^D,
\]
provided that $2A>2D+1+\delta$.

We now give a sufficient condition ensuring that the Euler products associated with realizations of our point process fail to belong to the Paley--Wiener space with positive probability.

\begin{theorem}\label{biquad}
Let $\alpha\leq1$, and suppose that the limit
\[
\lim_{R\to+\infty}R^\alpha\,\mathbb{E}_{\mathbb{P}}\,G_X(Rs)\,\overline{G_X(Rt)}=Q(s,t),
\] exists in the space $\mathscr{D}'(\mathbb{R}^2)$ of distributions on $\mathbb{R}^2$ and is nondegenerate in the sense that there exists a nonnegative infinitely differentiable function $\varphi$ with compact support on the positive half-line such that $Q(\varphi\otimes\varphi)>0$. Assume, furthermore, that the limit
\[
\lim_{R\to+\infty}R^{2\alpha}\,\mathbb{E}_{\mathbb{P}}\,G_X(Rs_1)G_X(Rs_2)\overline{G_X(Rt_1)}\,\overline{G_X(Rt_2)}=B(s_1,s_2,t_1,t_2).
\] exists in the space $\mathscr{D}'(\mathbb{R}^4)$ of distributions on $\mathbb{R}^4$. Then we have
\[
\mathbb{P}\bigl(\{X:G_X\notin L_2(\mathbb R)\}\bigr)>0.
\]
Under the additional assumption that the point process $\mathbb{P}$ has bounded deviation we also have
\[
\mathbb{P}\bigl(\{X:X\text{ is a uniqueness set for }PW\}\bigr)>0.
\]
\end{theorem}

Recall that two configurations are called symmetrically equivalent if there exists a compact subset of the phase space such that the two point sets coincide outside this compact set and have the same number of particles inside it. In the special case of the real line, two configurations are symmetrically equivalent if and only if a compactly supported diﬀeomorphism of the line takes one configuration to the other. Clearly, two symmetrically equivalent configurations on the line have the same excess (or deficiency) in the Paley--Wiener space. A Borel subset of the space of configurations that can be represented as a union of symmetric equivalence classes is called a symmetric event. The completion of the sigma-algebra of all symmetric events is called the symmetric sigma-algebra of our point process. Theorem \ref{biquad} now yields the following corollary.

\begin{corollary}
Suppose that the point process $\mathbb{P}$ has bounded deviation and satisfies the assumptions of Theorem \ref{biquad}. If the symmetric sigma-algebra of the point process $\mathbb{P}$ is trivial, then almost all its realizations are uniqueness sets in the Paley--Wiener space $PW$.
\end{corollary}

\begin{remark}
The fact that the excess of the sine process is equal to $1$ follows immediately from Theorem 1; a detailed proof will be given in the sequel to this note. Recently, Bourgade, Sodin, and Zeitouni announced an exact computation of the excess of the sine-$\beta$ process for all values of the parameter $\beta$. It would be very interesting to establish the Bourgade–Sodin–Zeitouni Theorem using the elementary methods developed in this note.
\end{remark}

Stochastic Euler products of point processes have been studied, in particular, in \cite{Assiotis}, \cite{bufetov-quasi}, \cite{CNN}, \cite{CHNNR}, \cite{NN}, and \cite{VV}.

\section{Proofs of Propositions \ref{small} and \ref{big}}

\subsection{Proof of Proposition \ref{small}}
First, recall that
\[
PW=\left\{f\in L_2(\mathbb R):\supp\hat f\subset[-\pi;\,\pi]\right\},
\]
and, by the Paley--Wiener Theorem in the form of Plancherel and Pólya, an entire function $f$ belongs to $PW$ if and only if
\begin{equation}\label{exp-gr-broad}
\limsup\limits_{r\to\infty}\frac{\log\max\limits_{z:|z|=r}|f(z)|}{r}\leq\pi
\end{equation}
and
\begin{equation}\label{l2r}
\int_{\mathbb R}|f(x)|^2\,dx<+\infty.
\end{equation}
\begin{remark}
For the purposes of this note, however, it suffices to use the original Paley--Wiener condition
\begin{equation}\label{exp-gr}
\sup_{z\in\mathbb C}\frac{|f(z)|}{\exp(\pi|z|)}<+\infty
\end{equation}
and its polynomially weighted extension \eqref{type-pi-poly}.
\end{remark}

We now verify the conditions of the Paley-Wiener Theorem for the Euler product associated with our configuration.

\textbf{1)} Let $r\in\mathbb N$ satisfy $r>B+N_0+19$. Choose the $r$ points $p_1,\ldots,p_r$ of our configuration $X$ that are closest to the origin, and set $Y=X\setminus\{p_1,\ldots,p_r\}$. The map $X\mapsto Y$ is Borel measurable. First observe that, under our assumptions, the Euler product $G_Y$ satisfies condition \eqref{exp-gr} almost surely with respect to $\mathbb{P}$. Indeed, by \eqref{rad-gr} and the Borel--Cantelli Lemma, for $\mathbb{P}$-almost every $X\in\mathscr{C}(\mathbb R)$ we have
\[
\limsup_{m\to+\infty}\frac{m^8}{e^{2\pi m}}\int_{-\pi}^{\pi}\left|G_Y(me^{i\theta})\right|^2\,d\theta<+\infty,
\]
which implies \eqref{exp-gr}.

\textbf{2)} Condition \eqref{l2r} for $G_Y$ follows directly from \eqref{sm-l2}.

\textbf{3)} It follows that, almost surely with respect to $\mathbb{P}$, the function $G_Y$ belongs to the Paley--Wiener space $PW=PW[-\pi;\,\pi]$.

\textbf{4)} If a function $G$ belongs to the Paley--Wiener space and
\[
\int_{\mathbb R}t^2|G(t)|^2\,dt<+\infty,
\]
then the function $t\mapsto tG(t)$ also belongs to the Paley--Wiener space. Our assumptions on $G_X$ imply that, for every set $q=\{q_1,\ldots,q_\ell\}$ of $\ell$ distinct particles of the configuration $X$, the function $G_{X^q}$ associated with the configuration $X^q=X\setminus\{q_1,\ldots,q_\ell\}$ is square-integrable on the real line:
\[
\int_{\mathbb R}|G_{X^q}(t)|^2\,dt<+\infty
\]
for $\mathbb{P}$-almost every configuration $X\in\mathscr{C}(\mathbb R)$. It follows that $G_{X^q}$ belongs to the Paley--Wiener class $PW=PW[-\pi;\,\pi]$, as required. \qed

We turn to the proof of Proposition \ref{big}.

\subsection{Proof of Proposition \ref{big}}
Let $N > A + 3$. From the Plancherel--Pólya Discrete Sampling Theorem \cite{PP37} it directly follows that if the point set $Y$ is obtained from $X$ by removing $N$ particles (in order to ensure Borel choice, remove $N$ positive particles closest to 0), then $G_Y \in PW$. Suppose that $\delta>0$, $|k|\geq r_0$, and
\[
\mathbb{P}\!\left(\left\{X:\left|\log|G_X(ik)|-\pi|k|\right|>A\log(2+|k|)\right\}\right)\leq\frac{C}{|k|^{1+\delta}}.
\]
Then, for $\mathbb{P}$-almost every $X$, we have
\[
\liminf\limits_{k\in\mathbb Z,\,|k|\to\infty}|G_X(ik)|\exp(-\pi|k|)\cdot(1+|k|)^{[A]+3}>0,
\]
and, by Proposition \ref{rtozinf}, also
\[
\inf\limits_{v\in\mathbb R,\,|v|\to\infty}|G_X(iv)|\exp(-\pi|v|)\cdot(1+|v|)^{[A]+3}>0.
\]
The Proposition now follows directly from Proposition \ref{low-imag-prop}. \qed

\subsection{Non-convergence in probability ensures completeness}

Let
\[
\begin{aligned}
\mathscr{E}_0(PW)=\bigl\{X\in\mathscr{C}(\mathbb R):\;&X\text{ is a uniqueness set for }PW\bigr\}.
\end{aligned}
\]

\begin{proposition}
Let $\mathbb{P}$ be a point process on $\mathbb R$ satisfying \eqref{rad-gr} and \eqref{log-dev}. Suppose, moreover, that there exist $a\in\mathbb R$, $\varepsilon>0$, and a sequence $t_n\to+\infty$ such that
\[
\liminf_{n\to\infty}\mathbb{P}\bigl(\{X:|G_X(t_n+ia)|>\varepsilon\}\bigr)>0.
\]
Then $\mathbb{P}(\mathscr{E}_0(PW))>0$.
\end{proposition}

\begin{proof}
By Proposition \ref{big}, the deviation of $X$ is almost surely finite. Suppose, to the contrary, that $\mathbb{P}(\mathscr{E}_0(PW))=0$. By Proposition \ref{bdd-dev}, it follows that $G_X\in PW$ almost surely. Consequently,
\[
\lim_{n\to\infty}G_X(t_n+ia)=0
\]
almost surely and, a fortiori, in probability. However, the convergence in probability to zero contradicts the assumption. The proof is complete.
\end{proof}

\section{Proof of Theorem \ref{biquad}}

\begin{proof}
Under the assumptions of the Theorem, there exists a nonnegative infinitely differentiable function $\varphi:\mathbb R\to[0,1]$ with compact support such that
\[
\int_{\mathbb R}\varphi(t)\,dt=1,\qquad \supp\varphi\subset[1,+\infty),
\]
and $Q(\varphi\otimes\varphi)>0$. Of course, we also have
\[
B(\varphi\otimes\varphi\otimes\varphi\otimes\varphi)\geq0.
\]
Our assumptions imply the existence of a positive limit
\[
\lim_{R\to+\infty}R^\alpha\,\mathbb{E}\left|\int G(Rs)\,\varphi(s)\,ds\right|^2=c_1>0,
\]
as well as the existence of the limit
\[
\lim_{R\to+\infty}R^{2\alpha}\,\mathbb{E}\left|\int G(Rs)\,\varphi(s)\,ds\right|^4=c_2.
\]
By the Cauchy--Bunyakovsky--Schwarz Inequality, we have $c_2\geq c_1^2>0$. The Paley–Zygmund Inequality now implies that there exist constants $c_3>0$ and $\rho_0>0$, depending only on $c_1$ and $c_2$, such that
\[
\liminf_{R\to+\infty}\mathbb{P}\!\left(\left\{X:\left|\int G_X(Rs)\,\varphi(s)\,ds\right|^2\geq c_3R^{-\alpha}\right\}\right)\geq\rho_0.
\]

Observing that
\[
\begin{aligned}
\left|\int_{\mathbb R}G_X(Rs)\,\varphi(s)\,ds\right|^2
&\leq\int_{\mathbb R}|G_X(Rs)|^2\,\varphi(s)\,ds\\
&=\frac1R\int_{\mathbb R}|G_X(u)|^2\,\varphi\!\left(\frac{u}{R}\right)du.
\end{aligned}
\]
we obtain
\[
\liminf_{R\to+\infty}\mathbb{P}\!\left(\left\{X:\int_{\mathbb R}|G_X(u)|^2\,\varphi\!\left(\frac{u}{R}\right)du\geq c_3R^{1-\alpha}\right\}\right)\geq\rho_0>0.
\]
Consequently, there exist $a,b\in\mathbb R$, $1<a<b$, and $c_4>0$ such that
\begin{equation}\label{key-est}
\liminf_{R\to+\infty}\mathbb{P}\!\left(\left\{X:\int_{aR}^{bR}|G_X(u)|^2\,du\geq c_4R^{1-\alpha}\right\}\right)\geq\rho_0>0,
\end{equation}
which directly implies, also in the case $\alpha=1$ (indeed, if a random function is square-integrable almost surely, then the integrals of its square over intervals far from zero converge in probability to zero, which contradicts \eqref{key-est}), that
\[
\mathbb{P}\bigl(\{X:G_X\notin L_2\}\bigr)>0.
\]
The proof is complete.
\end{proof}

\textbf{Acknowledgements.} I am deeply grateful to A.~S.~Arakcheev, A.~A.~Borichev, M.~A.~Borovikov, S.~M.~Gorbunov, M.~L.~Sodin, and I.~I.~Zavolokin for useful discussions.

The author is a winner of the ``Leader'' competition of the BASIS Foundation for the Advancement of Theoretical Physics and Mathematics, and would like to thank the jury and sponsors.

Part of this paper was written in the framework of a
research group at the Sirius Mathematical Centre. I am deeply grateful to the Sirius Mathematical Centre for its warm hospitality.


\begin{thebibliography}{99}

\bibitem{ahiezer}
N.~I.~Akhiezer, \emph{Lectures on the Theory of Approximation}, Gostekhizdat, 1947.

\bibitem{Akh54}
N.~I.~Akhiezer, On a generalization of the Fourier transform and the Paley--Wiener theorem. \emph{Doklady Akademii Nauk SSSR} \textbf{96} (1954), no.~4, 889--892.

\bibitem{Assiotis}
T.~Assiotis,
Random entire functions from random polynomials with real zeros,
\emph{Advances in Mathematics} \textbf{410} (2022), article 108701.
%\doi{10.1016/j.aim.2022.108701}; \arxiv{2202.03362}.

\bibitem{BB}
A.~Baranov, Yu.~Belov, A.~Borichev, \textit{Hereditary completeness for systems of exponentials and reproducing kernels}, Adv. Math. \textbf{235} (2013), 525--554.

\bibitem{bor-ok}
A.~Borodin, A.~Okounkov, A Fredholm determinant formula for Toeplitz determinants. Integr. Equ. Oper. Theory 37 (2000), 386--396.

\bibitem{bufetov-sine-fun}
Alexander Bufetov,
The Expectation of a Multiplicative Functional under the Sine-Process, Functional Analysis and Its Applications
Volume 58 (2024), 120--128.

\bibitem{BQS}
A.~Bufetov, Y.~Qiu, A.~Shamov, \textit{Kernels of conditional determinantal measures and the Lyons--Peres completeness conjecture},
arXiv:1612.0675, JEMS, 23(5), 2021, 1477--1519.

\bibitem{bufetov-excess}
A.~I.~Bufetov, \textit{The sine-process has excess one}, arXiv:1906.02098, 2019.

\bibitem{bufetov-quasi}
A.~I.~Bufetov, Quasi-symmetries of determinantal point processes. Ann. Probab. 46 (2018), 956--1003.

\bibitem{BS}
A.~I.~Bufetov and T.~Shirai,
Quasi-symmetries and rigidity for determinantal point processes associated with de Branges spaces,
\emph{Proceedings of the Japan Academy, Series A, Mathematical Sciences} \textbf{93} (2017), no.~1, 1--5.
%\doi{10.3792/pjaa.93.1}; \arxiv{1606.01564}.

\bibitem{CNN}
R.~Chhaibi, J.~Najnudel and A.~Nikeghbali,
\emph{A limiting random analytic function related to the CUE}, preprint, 2014.
arXiv:1403.7814.

\bibitem{CHNNR}
R.~Chhaibi, E.~Hovhannisyan, J.~Najnudel, A.~Nikeghbali and B.~Rodgers,
The limiting characteristic polynomial of classical random matrix ensembles,
\emph{Annales Henri Poincar\'e} \textbf{20} (2019), 1093--1119.
%\doi{10.1007/s00023-019-00769-4}; \arxiv{1707.09956}.

\bibitem{CNN-inv}
Reda Chhaibi, Joseph Najnudel, Ashkan Nikeghbali, ``The Circular Unitary Ensemble and the Riemann zeta function: the microscopic landscape and a new approach to ratios''.
Inventiones mathematicae 207 (2017), 23--113.

\bibitem{CL}
L.~Chistyakov, Yu.~Lyubarskii, \textit{Random perturbations of exponential Riesz bases in $L_2(-\pi,\pi)$}, Ann. Inst. Fourier \textbf{47} (1997), no.~1, 201--233.

\bibitem{geronimo}
J.~S.~Geronimo, K.~M.~Case, ``Scattering theory and polynomials orthogonal on the unit circle'', J. Math. Phys. 20:2 (1979), 299--310.

\bibitem{ghosh-exp}
S.~Ghosh, Determinantal processes and completeness of random exponentials: the critical case. Probab. Theory Relat. Fields 163 (2015), 643--665.

\bibitem{Gor25}
S.~M.~Gorbunov, Unitary transform diagonalizing the confluent hypergeometric kernel, \emph{Sb. Math.} \textbf{216} (2025), no.~12, 1643--1663.

\bibitem{johansson}
Kurt Johansson,
On random matrices from the compact classical groups.
Ann. of Math. (2), 145(3):519--545, 1997.

\bibitem{ibrag}
I.~A.~Ibragimov, On a theorem of G. Szeg\"o, Math. Notes, 3:6 (1968), 442--448.

\bibitem{kallenberg}
O.~Kallenberg, \emph{Random Measures}. Academic Press Inc., New York (1983).

\bibitem{KL38}
M.~I.~Kontorovich and N.~N.~Lebedev, On the method of solving certain problems of mathematical physics.
\emph{Journal of Experimental and Theoretical Physics (JETP)} \textbf{8} (1938), no.~10--11, 1192--1206.

\bibitem{Levin}
B.~Ya.~Levin, \textit{Distribution of Zeros of Entire Functions}, GITTL, Moscow, 1956.

\bibitem{lyons-det}
R.~Lyons, Determinantal probability measures. Publ. Math. Inst. Hautes \'Etudes Sci. (98) (2003), 167--212.

\bibitem{NN}
J.~Najnudel and A.~Nikeghbali,
\emph{Convergence of random holomorphic functions with real zeros and extensions of the stochastic zeta function},
preprint, 2022, revised 2024.
arXiv:2202.04284.

\bibitem{phrag}
E.~Phragm\'en (1904), ``Sur une extension d'un th\'eor\`eme classique de la th\'eorie des fonctions'', Acta Mathematica 28 (1904), 351--368.

\bibitem{phrag-linde}
L.~E.~Phragm\'en and E.~Lindel\"of (1908), ``Sur une extension d'un principe classique de l'analyse et sur quelques propri\'et\'es des fonctions monog\`enes dans le voisinage d'un point singulier'', Acta Math. 31 (1), 381--406.

\bibitem{pl-po}
M.~Plancherel, G.~P\'olya, Fonctions enti\`eres et int\'egrales multiples de Fourier. \emph{Commentarii Mathematici Helvetici} \textbf{9} (1936), 224--248.

\bibitem{PP37}
M.~Plancherel, G.~Pólya, Fonctions enti\`eres et int\'egrales de Fourier multiples. II.,\
\emph{Commentarii Mathematici Helvetici} \textbf{10} (1937), no.~2, 110--163.

\bibitem{polya}
G.~P\'olya, Untersuchungen \"uber L\"ucken und Singularit\"aten von Potenzreihen, Mathematische Zeitschrift 29 (1929), 549--640.

\bibitem{shirai}
T.~Shirai and Y.~Takahashi, Random point fields associated with certain Fredholm determinants I: fermion, Poisson and boson point processes. J. Funct. Anal. 205 (2003), 414--463.

\bibitem{sosh}
A.~Soshnikov, Determinantal random point fields. Russian Math. Surveys (Uspekhi Mat. Nauk 55:5, 1071--1076) 55(5) (2000), 923--975.

\bibitem{VV}
B.~Valk\'o and B.~Vir\'ag,
The many faces of the stochastic zeta function,
\emph{Geometric and Functional Analysis} \textbf{32} (2022), 1160--1231.

\end{thebibliography}
\end{document}